\documentclass{amsart}
\newif\ifpdf\ifx\pdfoutput\undefined\pdffalse\else\pdftrue\fi
\ifpdf\usepackage[colorlinks,pagebackref]{hyperref} 
\else\usepackage[hypertex,pagebackref]{hyperref}\fi 
\usepackage{graphicx}
\usepackage{amscd}
\usepackage[curve]{xypic}

\newcommand{\ocomment}[1]{}
\newbox\mybox
\def\overtag#1#2#3{\setbox\mybox\hbox{$#1$}\hbox to
  0pt{\vbox to 0pt{\vglue-#3\vglue-\ht\mybox\hbox to \wd\mybox
      {\hss$\ss#2$\hss}\vss}\hss}\box\mybox}
\def\undertag#1#2#3{\setbox\mybox\hbox{$#1$}\hbox to 0pt{\vbox to
    0pt{\vglue#3\vglue\ht\mybox\hbox to \wd\mybox
      {\hss$\ss#2$\hss}\vss}\hss}\box\mybox}
\def\lefttag#1#2#3{\hbox to 0pt{\vbox to 0pt{\vss\hbox to
      0pt{\hss$\ss#2$\hskip#3}\vss}}#1}
\def\righttag#1#2#3{\hbox to 0pt{\vbox to 0pt{\vss\hbox to
      0pt{\hskip#3$\ss#2$\hss}\vss}}#1}
\let\ss\scriptstyle

\def\Dot{\lower.2pc\hbox to 2.5pt{\hss$\bullet$\hss}}
\def\Circ{\lower.2pc\hbox to 2.5pt{\hss$\circ$\hss}}
\def\Vdots{\raise5pt\hbox{$\vdots$}}
\def\splicediag#1#2{\xymatrix@R=#1pt@C=#2pt@M=0pt@W=0pt@H=0pt}

\renewcommand\frame[2][3pt]{\hbox{$\vcenter{\hbox{\vrule\vbox {\hrule\kern#1%
  \hbox{\kern#1$#2$\kern#1}\kern#1\hrule}\vrule}}$}}

\newcommand\lineto{\ar@{-}}
\newcommand\dashto{\ar@{--}}
\newcommand\dotto{\ar@{.}}
\newcommand{\C}{{\mathbb C}}
\newcommand{\Z}{{\mathbb Z}}

\newtheorem{theorem}{Theorem}[section]
\newtheorem*{theorem*}{Theorem}
\newtheorem*{ecconjecture*}{End-Curves Conjecture}
\newtheorem*{ciconjecture*}{Casson Invariant Conjecture}

\newtheorem{proposition}[theorem]{Proposition}

\theoremstyle{definition}

\newtheorem{example}[theorem]{Example}

\begin{document}
\title[Topology, geometry, and equations of surface singularities]
{Topology, Geometry, and Equations of Normal Surface Singularities}
\author{Jonathan Wahl}
\thanks{Research supported under NSA grant no.\ H98230-04-1-0053}
\address{Department of Mathematics\\The University of
North Carolina\\Chapel Hill, NC 27599-3250}
\email{jmwahl@email.unc.edu}
\dedicatory{To Gert-Martin Greuel on his sixtieth}
\begin{abstract}

  In continuing joint work with Walter Neumann, we consider the relationship between
  three different points of view in
describing a (germ of a) complex normal surface singularity. The
explicit equations of a singularity allow one to talk
about hypersurfaces, complete intersections, weighted homogeneity,
Hilbert function, etc. The geometry of the singularity could involve analytic
aspects of a good resolution, or existence and properties of Milnor fibres;
one speaks of geometric genus, Milnor number, rational singularities,
the Gorenstein and $\mathbb Q$-Gorenstein properties, etc.  The topology of the
singularity means the description of its link, or equivalently (by a
theorem of Neumann) the configuration of the exceptional curves
in a resolution.
   We survey ongoing work
    (\cite{nw1},\cite{nw2}) with
    Neumann to study the possible geometry and equations when
    the topology of the link is particularly simple, i.e.
    the link has no rational homology, or
    equivalently the exceptional configuration in a
    resolution is a tree of rational curves.  Given such a link, we
    ask whether there exist
    ``nice'' singularities with this topology.  In our situation, that
    would ask if
    the singularity is a quotient of a special kind of explicitly
    given complete
    intersection (said to be ``of splice type'') by an explicitly
    given
    abelian group; on the topological level, this quotient gives the universal
    abelian cover of the link.   Our major result gives a
    topological condition (i.e., a condition on the resolution graph)
    that there exists a singularity which arises in this way (and hence one whose equations can be
    written ``explicitly'').  T.
    Okuma (\cite{okuma2}) has recently proved our Conjecture
    that
    rational and minimally
    elliptic singularities are all ``splice-quotients''.   We summarize
    first the well-studied case of plane curve
    singularities, to see what one might mean about geometry,
    topology, and equations in that case.  There follows an
    introductory discussion of normal surface singularities, before considering our recent work.

    The purpose of the article is to survey the main ideas and
    directions, rather than to describe details, which can be found
    in other papers such as \cite {nw1}.

\end{abstract}
\maketitle

\section{Introduction}
   To understand what we mean by ``topology, geometry, and
   equations,'' we start with the germ at the origin of a complex irreducible
   (and
   reduced) plane curve singularity
   $C=\{f(x,y)=0\}\subset \mathbb C ^{2} $. Intersecting with a small
   $3$-sphere gives a knot $L$ in the 3-sphere.  The embedded topology
   of the knot was studied by K. Brauner and E. K\"{a}hler in the
   1920's (see \cite{n} for some details).  Their approach
   resulted in the topological description of the knot by iterated
   cabling on a torus knot.  The description is given by a sequence
   of pairs of positive integers $(p_{i},r_{i})$, the Puiseux pairs, which
   can be read off
   a fractional power series which parametrizes the curve;
   equivalently, a related approach produces the sequence of Newton pairs
   $(p_{i},q_{i})$.  A point is that the ``link'' of the singularity
   is intrinsically just a circle, so the ``topology'' of the
   situation should mean the embedded topology (i.e., knot type of
   $L$).  In this case, from
   the equation $f(x,y)=0$ one can iteratively read off the Puiseux pairs
    using Newton diagrams, and this data describes
   fully the topology.
   From a more geometric point of view, consider an embedded resolution of
   the curve by
     blowing up $\mathbb C^{2}$, until the reduced total transform of
     the curve has normal crossings.  Again, the numerical data described
     above can be
     read off from the configuration of the exceptional curves
     and their self-intersections, plus the intersection with the transform
     of $C$.
       Alternatively, one may consider the \emph{value semigroup}
     of the singularity; writing the integral closure of the local
     ring of $C$ as $\mathbb C\{\{t\}\}$, one considers the collection
     of $t$-orders of all
     elements of the subring.  Then this value semigroup is equivalent
     to the data of the Puiseux pairs.

     What about recovering the equation of the curve from the above
     data?  Of course, one has \emph{equisingular}
   families for which the embedded
     topology is constant (i.e., same numerical data), but with
     analytically distinct individual curves.  So one would ideally like to
     write down every plane curve singularity with given topological
     type.  There are several ways to do this.  For instance, in
     the Appendix to Zariski's book
     (\cite{zariski}), B.
     Teissier
     considers a monomial curve given by generators of the value
     semigroup; this is known to be a complete intersection, and is
     weighted homogeneous.  The versal deformation is smooth
     and also carries a $\mathbb C^*$-action; then the deformations of
     non-negative weight give all curves with the same value
     semigroup (though many of these curves are now no longer
     planar).  Another approach is to write down the most general Puiseux
     series with given Puiseux pairs, as in \cite{zariski} Chapitre III (cf.
     also \cite{eisenbud-neumann}, Appendix to Chapter 1).
     For one
     Puiseux pair $(p,q)$ a family containing every analytic type is given by:
     $$x^{p}+y^{q}+\sum t_{ij}x^{i}y^{j}=0,$$
     the sum over $(i,j)$ such that $$i/p\ +\ j/q\ >1,\ \  0<i<p-1,\ 0<j<q-1.$$

     For a reducible plane curve, one must keep track
     not only of the topological type of each branch, but the linking
     numbers as well.

     A new and important
     development in the study of curves appeared in the use of the
     \emph{splice diagrams} of Siebenmann, by D. Eisenbud and W. Neumann
     \cite{eisenbud-neumann}.  These will be discussed below to study surface
     singularities.  For this and other topics in the topology of plane
     curve singularities, we refer again to \cite{n}.

     The moral for us is that at least for irreducible plane curve
     singularities, we know how to recover all the relevant (embedded)
     topological data from the equation, or a resolution, or the value
     semigroup; we can tell when a set of data comes from a singularity;
     and from the data
      we can write down equations of all plane
     curves of that topological type.

     \section{The basics of normal surface singularities}

     Working up a dimension from the case of curves, consider now
     the germ at the origin of an
     an isolated hypersurface singularity  $Y=\{f(x,y,z)=0\}\subset
     \mathbb C^3$.  The \emph{link} $\Sigma$ of the singularity
     is the intersection with a
     small $5$-sphere centered at the origin.  $\Sigma$ is a compact
     connected oriented 3-manifold,
     knotted somehow in $S^5$.  The local topology of the pair
     $(Y,\mathbb C^3)$ is given by the topological cone over
     $(\Sigma, S^5)$; in particular, $Y$ is a topological manifold
     at the origin iff $\Sigma$ is homeomorphic to the $3$-sphere.
     As for knotting in the 5-sphere, one has (via the map $f/|f|$)
      the \emph{Milnor fibration}; the complement of $\Sigma$
      in $S^5$ fibres over the circle, and the fibres are
      $4$-manifolds with boundary $\Sigma$ and with homology only in
      the middle dimension.  The rank of this second
      homology group, called the \emph{Milnor number $\mu$}, is
      known to be computable as the colength of the Jacobian ideal
      $(f_x,f_y,f_z)$ in $\mathbb C[x,y,z]$.  This story is
      described in the classical book of Milnor
      \cite{milnor-book}.

      The whole subject of the topology of normal surface
      singularities really began with an important discovery by D. Mumford
      in 1960.  He showed that if the link $\Sigma$ is
      simply-connected, then not only is it the $3$-sphere (as the
      Poincar\'{e} Conjecture asserts), but it is unknotted in
      $S^5$, and in fact the origin is a non-singular point.
      Mumford's argument works not only in the hypersurface case.
      Suppose one has a germ of a normal surface singularity (from now
      on: NSS)
      $(Y,0)\subset(\mathbb C^n,0)$, and one considers the
      intersection with a small sphere $\Sigma=Y\cap
      S^{2n-1}$.  The theorem asserts that if $\Sigma$ is simply-connected,
      then $Y$ is smooth at $0$.
      Therefore, unlike in the case of curves,  for a NSS
       the topology of the link can give you a huge amount of
       information about the singularity; and in the simply
       connected case, it tells you \emph{everything} about the
      geometry of the point.  So, from now on, by \emph{the topology}
      of a NSS we shall mean simply the topological type of the 3-manifold
      $\Sigma$.


   The natural way to see the topology of a NSS $(Y,0)$ is via
    a ``good" resolution $\pi:(\tilde Y ,E)\rightarrow (Y,0)$.
   Thus, $\tilde Y$ is smooth, $\pi$ is proper and maps $\tilde Y -E$
   isomorphically
   onto $Y-\{0\}$, and
   $E=\pi^{-1}(0)$ is a divisor consisting of smooth projective curves $E_i$,
   intersecting transversally (there is in fact a \emph{minimal} good
   resolution, in an obvious sense).  One can associate to $E$ in
   the usual way the weighted resolution dual graph $\Gamma$: each irreducible component
   $E_{i}$ of $E$
   gives a vertex, intersection points give edges of the graph, and each
   vertex is weighted
    by the degree of
   the normal bundle of the corresponding irreducible curve.  In
   addition, the graph $\Gamma$ is decorated at each vertex with the genus
   of the corresponding curve.  The link $\Sigma$ can be
   reconstructed from $\Gamma$ (it is a \emph{graph manifold});
   this is because $\Sigma$ may be viewed as the boundary of a
   tubular neighborhood of $E$ on the smooth surface $\tilde Y$.  A critical
   fact, noted originally by
   P. DuVal, is that the intersection matrix $(E_i \cdot E_j )$ is
   negative-definite.  We
   conclude that the first betti number of $\Sigma$ equals the number of
   cycles in the graph plus twice the sum of
   the genera of all the $E_i$.

   If $\Sigma$ is simply connected, then
   \emph{a fortiori} $H_1(\Sigma;\mathbb Q)=0$,
   i.e., $\Sigma$
   is a \emph{rational homology sphere} (which we denote by $\mathbb
   Q$HS.)
    Thus, $E$ is a tree of smooth rational
   curves.
    Mumford shows in this case how to compute
   the fundamental group in terms of
   the loops surrounding the exceptional components in $\tilde Y$.  One must
   use simple connectivity
   to show that $E$ can be
   contracted to a smooth point.

   Mumford's method soon led to some generalizations.  If
   $\pi _1 (\Sigma)$ is finite, then the singularity itself is a
   quotient $\mathbb C ^2 /G$, where $G\subset GL(2,\mathbb C)$ is
   a finite subgroup containing no pseudo-reflections.  It also
   turns out that when $\pi _1 (\Sigma)$ is solvable, one is in a
   very rigid (and well-understood) situation.  But in analogy with
   the case of curves (and thus, e.g., of singularities defined by
   $z^{n}=f(x,y)$), one expects that rarely will $\Sigma$ by itself
   determine the full analytic type of a singularity.  In fact, H.
   Laufer (\cite{laufert}), following some earlier results of G. Tjurina,
   gives a complete list of resolution dual graphs which have a
   unique analytic representative
   (he calls such singularities \emph{taut}); they are all rational or minimally
   elliptic (see below for definitions).

One should mention that work of H. Grauert shows that every
negative-definite weighted dual graph (with the genera included)
does indeed arise from resolving some NSS $(Y,0)$.  One pastes
together analytically a smooth surface with the desired curve
configuration, and proves that you can blow-down the curve
configuration to a point, necessarily on a normal analytic
surface.  (A general result of Hironaka proves that isolated
analytic singularities are algebraic.)  Of course, this purely
existential result gives no indication whether such a singularity
could have a nice property, like being a hypersurface or complete
intersection (other than an obvious numerical condition that
allows for a ``canonical divisor'' to exist--see below).  An
important result of W.Neumann \cite{neumann81} also shows that the
homeomorphism type of $\Sigma$ uniquely determines the graph
$\Gamma$ of the minimal good resolution (with the well-known
exceptions of cyclic quotient singularities and ``cusp"
singularities, where orientation must be taken into account).  The
bottom line is that the topology of NSS's is reflected exactly by
the dual graphs $\Gamma$.

 So, a natural question is what statements can be made about the
 (many) singularities with given topology.  One would not mind if the
 situation were similar to that of equisingularity for
plane curves, i.e., if all such analytic types
fit into nice topologically trivial (and
geometrically similar) families.
Unfortunately, that is not the case in general.  For instance,
recall that a germ $(Y,0)$ is called
 \emph{Gorenstein} if there exists a nowhere-0 holomorphic 2-form
 $\omega$ on $Y-\{0\}$ (and hence on $\tilde{Y}-E$).  (There is also
 the traditional
 purely algebraic definition in terms of the local ring.)  Such an $\omega$
 extends meromorphically over $E$, hence gives rise to a ``canonical
 divisor''$K$, an integral combination of the exceptional curves.  Complete
 intersections are Gorenstein, as are quotient singularities $\mathbb
 C^{2}/G$ as above if and only if $G\subset SL(2,\mathbb C)$.  For a
 dual graph $\Gamma$ to come from a Gorenstein singularity, it is
 necessary that there exist
 an integral divisor $K$ satisfying the adjunction
 rules $$K\cdot E_{i} +E_{i}\cdot E_{i}=2g(E_{i})-2$$
 for all $i$; we can call such a graph \emph{numerically Gorenstein}.
 A key result of Laufer \cite{laufermin} shows that for any graph which does
 not
 correspond to a rational double point or minimally elliptic singularity,
there always exists a non-Gorenstein singularity with that graph.
So, if you consider the link of \emph{any} hypersurface
singularity in $\mathbb C ^3$ which is not rational or minimally
elliptic (e.g., has multiplicity at least 4), then there exists a
non-Gorenstein singularity with the same link.  It is completely
unknown if a numerically Gorenstein graph always arises from a
Gorenstein singularity, even when $\Gamma$ is a rational tree.

Now, some non-Gorenstein singularities are still quite pleasant.
We call a singularity $(Y,0)$
\emph{$\mathbb Q$-Gorenstein} if the canonical line bundle on
$Y-\{0\}$ has finite order, i.e. some $r$-th tensor power has a
nowhere-0 section.  By taking cyclic covers, one sees equivalently
that $Y$ is a quotient of a Gorenstein
singularity by a finite abelian group.  For instance, all rational
singularities are
$\mathbb Q$-Gorenstein.  One fully expects, as a generalization to
Laufer's result above, that for any graph
which is not rational or minimally elliptic, there exist
non-$\mathbb Q$ -Gorenstein representatives.  (However, we have not
seen a proof of this statement.)

Finally, to illustrate how very
different singularities can have the same link, consider the
Brieskorn singularities
$$V(p,q,r)=\{x^p +y^q +z^r =0\}\subset \mathbb C^3,$$ and their
links $\Sigma(p,q,r)$.  Then $\Sigma(3,4,12)$ is homeomorphic to
$\Sigma(2,7,14)$ (both have as dual graph a single curve of genus
$3$, with weight $-1$).  Yet these singularities have different
multiplicities and different Milnor numbers ($66$ and $78$,
respectively).  Further, if $C$ is a smooth projective genus 3
curve, and $P\in C$, then
$$Y=\text{Spec} \bigoplus \Gamma (\mathcal{O}_{C}(nP))$$ has the same link as
these last 2, but
is $\mathbb Q$-Gorenstein  iff the degree 0 line bundle
$\mathcal{O}(K-4P)$ has finite order.

\section{Some cases when topology implies nice geometry}
By the ``geometry'' of a singularity $(Y,0)$, one is interested in
analytic issues which go beyond the ``mere'' topology of the link.
 Relevant notions include: embedding dimension and Hilbert function; complete
intersection, or Gorenstein, or $\mathbb Q$-Gorenstein; the
geometric genus $p_{g}=\text{dim}\ R^{1}\pi _{*}\mathcal O
_{\tilde{Y}}$, where $\pi :\tilde{Y}\rightarrow Y$ is a
resolution; nature of the defining equations and their syzygies.
Complete intersection singularities have a simply connected Milnor
fibre, and hence a Milnor number $\mu$. By an old formula of
Laufer \cite{laufermu}, in the complete intersection case one has
a relation between these invariants: $\mu -12p_g$ is an explicit
purely topological invariant.  As already indicated, topologically
equivalent germs could have very different geometry.  But suppose
one has a family $\mathcal Y \rightarrow T$, which has a
simultaneous equitopological resolution $\tilde{\mathcal Y}
\rightarrow \mathcal Y \rightarrow T$ (that is, one has a locally
trivial deformation of the exceptional sets). Then the geometric
genus is known to be constant; and quite generally, a small
deformation of a complete intersection (or Gorenstein) singularity
has the same property.  On the other hand, in analogy with what is
known from deforming space curves, one can easily have a
topologically trivial and simultaneous resolution family of
complete intersection surface singularities for which the
embedding dimension jumps (\cite{nemethi et al} gives nice
examples). Thus, in the general theory we
 should not be worried about jumping multiplicity or embedding
 dimension in ``geometrically nice'' families.

We review the two situations where one knows a great deal about
the
 geometry of
 a singularity $(Y,0)$ just from its topology (i.e., graph
 $\Gamma$).  A good general reference for the following discussion
 is \cite{reid}, Chapter 4.
 Consider the minimal good resolution
 $(\tilde{Y},E)\rightarrow (Y,0)$, and exceptional cycles
 $Z=\Sigma n_i E_i$.  Among all such divisors, there is a minimal non-0 cycle $Z_{o}$
 with the property that $Z_{o}\cdot E_i \leq 0$, all $i$; this is the
 \emph{fundamental cycle} (also known as the \emph{numerical cycle} in
\cite{reid}), which is easily  (but not necessarily
quickly) computed.  The canonical line bundle $K$ of $\tilde{Y}$
satisfies $K\cdot E_i +E_i\cdot E_i
=2g(E_i)-2$ for every exceptional curve, hence one can
make sense of $K$ dotted with any cycle.  In
particular, $Z_{o}$ and $K$ are computed from the graph.

We say that $Y$ is a \emph{rational singularity} if $Z_{o}\cdot
(Z_{o}+K)=-2$; a basic theorem says that this condition is
equivalent to the vanishing of the geometric genus (which is a
priori an analytic as opposed to topological invariant).  Seminal
work of M. Artin and E. Brieskorn, expanded upon by J. Lipman,
shows the multiplicity of a rational $Y$ is $m\equiv -Z_{o}\cdot
Z_{o}$, and the embedding dimension is $m+1$ (so the Hilbert
function is $H(n)=mn+1$).  $Y$ can be minimally resolved simply by
a sequence of blowing-up the singular points (which are always
themselves rational), and the exceptional divisor is a tree of
smooth rational curves.  The local Picard group (=divisor class
group) has finite order.  So the canonical line bundle on
$Y-\{0\}$ has finite order, from which it follows that $Y$ is
$\mathbb Q$-Gorenstein (though it is Gorenstein only for the
rational double points).  It was proved in \cite{wahl} that $Y$ is
defined by quadratic equations, and all the higher syzygies are
``linear'' in an appropriate sense.  All finite quotient
singularities $\mathbb C ^{2}/G$ are rational; writing down the defining equations is
a calculation in invariant theory.  Up to now, given the
graph of a rational singularity, there is no general method for
writing down explicit equations of a corresponding rational
singularity.  The best results, due to De Jong and van Straten
\cite{ds}, show how to do this if the rational graph has
reduced fundamental cycle.

Next, we say $Y$ is  \emph{minimally elliptic} if on the minimal
resolution  $Z_{o}\equiv -K$.  These singularities were introduced
by H. Laufer in \cite{laufermin}, though many of the results were
discovered independently by Miles Reid (in unpublished notes, but
see \cite{reid} ).  The definition is equivalent to $Y$ being
Gorenstein of geometric genus 1.  Except for cones over elliptic
curves and ``cusp'' singularities (whose resolution dual graph is
a cycle of smooth rational curves), the minimal good resolution
graph is a rational tree. When $m\equiv -Z_{o}\cdot Z_{o}$ is 1,2,
or 3, one has a hypersurface of multiplicity 2, 2, or 3
respectively in $\mathbb C ^{3}$; when $m\geq 4$, $m$ is the
multiplicity, the Hilbert function is $H(n)=mn$.  Further, $Y$ is
defined by quadratic equations, with ``linear syzygies'' except at
the last step \cite{wahl}.  Up to now, there is no general method
to write down equations, given the minimally elliptic resolution
graph.

We mentioned above that Laufer's result, and a presumed
generalization, would imply that \emph{any} resolution graph which
is neither a rational nor minimally elliptic graph would arise
from at least one non-$\mathbb Q$-Gorenstein singularity.  In
other words, we should be completely finished with the ``nice''
cases where the topology automatically forces some basic facts
about the geometry.

But there is one more situation in which a great deal can be said
from the topology---that is, if one additionally knows that
$(Y,0)$ is weighted homogeneous (that is, quasi-homogeneous, or
admits a good $\mathbb C^{*}$-action).  Then $Y=\text{Spec}\ A$,
where $A$ is a positively graded $\mathbb C$-algebra.  It follows
from early work of Orlik-Wagreich \cite{O-W} that (except for
cyclic quotient singularities) the exceptional divisor on the
minimal good resolution consists of one central smooth curve (Proj
$A$), and chains of smooth rational curves emanating from at least
3 points of this curve.  Put another way, consider the weight
filtration $\{I_{n}\}$ of $A$, where $I_{n}$ is the ideal
generated by elements of weight $\geq n$.  Take the weighted
blow-up $Z=\text{Proj}\  \oplus I_{n} \rightarrow Y=\text {Spec}\
A$ (the so-called \emph{Seifert partial resolution}).  Then $Z$ is
a normal surface with several cyclic quotient singularities along
its exceptional divisor, which is isomorphic to Proj $A$.  In
particular, $(Y,0)$ determines the following data: the isomorphism
class of the central curve; its conormal line bundle; the location
of the points on the curve at which $Z$ has a singularity; and the
data of the cyclic quotient singularities at these points.
Conversely, it was shown by H. Pinkham \cite{pinkham} and
independently by I. Dolgachev how to write down explicitly the
graded algebra $A$ from this data.  In other words, this data
uniquely determines the analytic type of the singularity.

Now suppose we have a weighted homogeneous singularity $Y$ with
\emph{rational} central curve.  Then the data you need to write
down $Y$ is numerical, contained in the graph $\Gamma$, except for
the (analytically significant) location of the intersection points
on the central rational curve.  Such singularities all have
$\mathbb Q$HS links, are rarely rational singularities.  For
instance, any Brieskorn hypersurface $V(p,q,r)$ for which $p$ is
relatively prime to $qr$ give such examples; but if $p,q,r \geq
4$, then these are neither rational nor minimally elliptic.  On
the other hand, due to the grading, it is not too hard to find the
numerical data that tell you when $Y$ is Gorenstein; and it turns
out that such $Y$ are always $\mathbb Q$-Gorenstein.  But much
more is true; we need to explain first some general facts.

\section{Universal abelian covers of singularities with $\mathbb Q
HS$ links}

There is an alternate way to describe the equations
of a weighted homogeneous singularity whose link is a $\mathbb
Q$HS.  Note that in general if the link $\Sigma$ of a singularity
is a $\mathbb Q$HS, then the first homology $H_{1}(\Sigma;\Z)$ is a
finite group computed directly from $\Gamma$; this \emph{discriminant
group} $D(\Gamma)$ is the cokernel of
the intersection pairing $(E_i \cdot E_j)$ (so, the order is the
absolute value of the determinant).  The universal abelian
covering $\tilde \Sigma \rightarrow \Sigma$ is finite, and it is
important to note that it can be realized by a finite map of germs
of NSS's $(X,0)\rightarrow (Y,0)$, unramified off the singular
points. We abuse notation and refer to the map $X\rightarrow Y$ as
the \emph{universal abelian covering of the singularity}, or the UAC.

Recall that a Brieskorn complete intersection (or BCI) is a
singularity $V(p_{1},\ldots,p_{t})$ defined by
$$\sum_{j=1}^{t}a_{ij}z_{j}^{p_{j}} =0,\
i=1,\dots, t-2\, ,$$ where $p_{i}\geq 2$ and every maximal minor
of the matrix $(a_{ij})$ has full rank.

\begin{theorem}\label{UAC} \cite{neumann83} Let $(Y,0)$ be a weighted homogeneous
singularity whose link is a $\mathbb Q HS$.  Then the UAC of
$(X,0)$ is a BCI as above.  The $p_{i}$ and the diagonal action of
the discriminant group on the ambient variables are explicitly
computed from the graph $\Gamma$ of $Y$.
\end{theorem}

Of course, the values of the $a_{ij}$  depend on the analytic
class of $t$ points on the central $\mathbb P ^{1}$.  Note that
the Theorem allows one to write down  ``explicit'' equations for
the singularity $(Y,0)$.  First, write down monomials generating
the ring of invariants for the action of the discriminant group
acting on $\mathbb C [z_1,\cdots,z_t]$; then, mod out by relations
on these monomials; finally, divide by relations which follow from
the BCI equations.  While not as direct a method for writing
equations as the aforementioned Pinkham-Dolgachev approach, this
result actually generalizes to many other cases.

\begin{example}\label{easy}
 The $E_{7}$ singularity (a rational double
point) has graph $\Gamma$ which is the Dynkin diagram for $E_{7}$.
Its discriminant group has order 2, and the UAC is $V(2,3,4)$,
which is the $E_{6}$ singularity; the group action on
$x^{2}+y^{3}+z^{4}=0$ is given by $(x,y,z)\mapsto (-x,y,-z)$.  So
the quotient is generated by the invariants
$A=x^{2},B=xz,C=z^{2},D=y$, with equations $AC-B^{2}=0$ from the
group action, and $A+D^{3}+C^{2}=0$ from the equation. This yields
the familiar (non-Brieskorn) equation for $E_{7}$:
$$B^{2}+C(C^{2}+D^{3})=0.$$
\end{example}

There are two striking aspects of Theorem \ref{UAC}.  First, the UAC
turns out to be not only Gorenstein (which is clear), but
even a complete intersection.  Second, it is a complete
intersection of very special type---a BCI (which need \emph{not} have $\mathbb Q HS$ link).

For many years we wondered what is the most general natural context
in which to view the previous theorem.  It thus became
clear one should study the UAC of singularities $(Y,0)$ with $\mathbb Q HS$
link, but whose graph $\Gamma$ has more than one node (unlike in the
weighted homogenous case).  If the UAC were to be a complete intersection, then $Y$
must be $\mathbb Q$-Gorenstein.  But
what generalization of BCI would work as UAC's for a wide class of $Y$'s, such
as rational singularities?

One considers first the case that the link is a $\mathbb Z HS$,
since then the UAC of the singularity is itself.  This was the
topic of several papers with Neumann \cite {neumann-wahl90, nw2}.
The discovery of the Casson invariant $\lambda (\Sigma)$ of a
three-dimensional $\mathbb Z HS$ around 1986 and its subsequent
calculation for certain examples led to the Conjecture of
Neumann-Wahl:

\begin{ciconjecture*}\cite{neumann-wahl90} For a complete intersection singularity with
 $\mathbb Z HS$ link,
the Casson invariant is one-eighth the signature of the Milnor
fibre.
\end{ciconjecture*}

Since the Casson invariant is topological, such a result would
imply that for a complete intersection (a very strong geometric
property), the $\mathbb Z HS$ topology determines the signature of
the Milnor fibre (and hence, by well-known formulae, also the
Milnor number and geometric genus).  Implicit in the Conjecture
(and what makes it provocative) is the prediction that the Milnor
fibre itself is somehow canonically associated to the link
(perhaps with some extra structure).

 In spite of some progress on this Conjecture \cite{nw2},
and counterexamples showing it does not generalize to links of
Gorenstein singularities (which might not even be smoothable)
\cite{nemethi et al}, the question remains open. However, it did
raise the problem of trying to write down explicit examples of
complete intersection singularities with $\mathbb Z HS$ links,
beyond the BCI's $V(p_{1},\ldots,p_{t})$ (where the $p_{i}$ are
pairwise relatively prime).  This ultimately led to the discovery
of \emph{complete intersections of splice type}, which play a role
in the general problem.

\section{splice diagrams and complete intersections of splice
type--$\mathbb Z HS$ case}

The usual topological description of a singularity link is via
plumbing according to the resolution dual graph $\Gamma$.  But
when the link is a $\mathbb Z HS$ (i.e., the intersection matrix
is unimodular), there is another topological construction, from a
different kind of graph, which can be computed from $\Gamma$.  The
following discussion is largely taken from \cite{nw2}, itself
depending heavily upon \cite{eisenbud-neumann}.

Suppose first that $K_{i}$ is a knot in a $\mathbb Z HS$  $\Sigma_{i},
i=1,2$. Then one may ``splice'' the two three-manifolds together
along the knots to form a new
$\mathbb Z HS$:  remove from each $\Sigma _{i}$ a tubular
neighborhood of $K_{i}$, and then paste together along the
boundaries (which are tori), but switching the roles of meridian and
longitude.  (Of course, orientation needs to be handled carefully).

A \emph{splice diagram} is a finite tree with vertices
only of valency 1 (``\emph{leaves}'') or $\ge3$ (``\emph{nodes}'') and
with a collection of integer weights at each node, associated to the
edges departing the node. The following is an example:
$$\splicediag{12}{30}{
\Circ&&&\Circ\\
&\Circ\lineto[ul]_(.25){2}\lineto[dl]^(.25)3
&\Circ\lineto[dr]_(.25){5}\lineto[ur]^(.25){2}
\lineto[l]_(.2){11}_(.8){7}\\
\Circ&&&\Circ
}$$
For an edge connecting two nodes in a splice diagram the \emph{edge
  determinant} is the product of the two weights on the edge minus the
product of the weights adjacent to the edge. Thus, in the above
example, the one edge connecting two nodes has edge determinant
$77-60=17$.  This example is supposed to represent the result of
splicing together the Brieskorn homology spheres $\Sigma(2,3,7)$
and $\Sigma(2,5,11)$ along the knots obtained by setting the last
coordinate equal to 0 in the defining equations.  Each leaf of a
splice diagram corresponds to a knot on the corresponding $\mathbb Z
HS$.

The splice diagrams that classify integral homology sphere singularity
links satisfy the following conditions on their weights:
\begin{itemize}
\item the weights around a node are positive and pairwise coprime;
\item the weight on an edge ending in a leaf is $>1$;
\item all edge determinants are positive.
\end{itemize}


\begin{theorem}[\cite{eisenbud-neumann}]
The integral homology spheres that are singularity links are in one-one
correspondence with splice diagrams satisfying the above conditions.
\end{theorem}

The splice diagram and resolution diagram for the singularity
determine each other uniquely, and indicate how to construct the
link by splicing or by plumbing.  To go from resolution to splice
diagram, one collapses all vertices of valency 2, and uses as
weights the absolute value of the intersection matrices of certain
subdiagrams. (Computing in the other direction is harder, and
given in \cite{eisenbud-neumann} or an appendix to \cite{nw2}).
The example above corresponds to the resolution diagram
$$
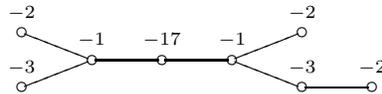
\xymatrix@R=6pt@C=24pt@M=0pt@W=0pt@H=0pt{
\\
\overtag{\Circ}{-2}{8pt}&&&&\overtag{\Circ}{-2}{8pt}\\
&\overtag{\Circ}{-1}{8pt}\lineto[ul]\lineto[dl]\lineto[r]&
\overtag{\Circ}{-17}{8pt}&\overtag{\Circ}{-1}{8pt}\lineto[ur]\lineto[dr]\lineto[l]&\\
\overtag{\Circ}{-3}{8pt}&&&&\overtag{\Circ}{-3}{8pt}\lineto[r]&\overtag{\Circ}{-2}{8pt}}
$$
An important point is that the ends of the diagrams (in this
$\mathbb Z HS$ case) correspond to certain natural isotopy classes
of knots in the 3-manifold.

A surprising early discovery of the Neumann-Wahl collaboration was that if a
singularity link is a $\mathbb Z HS$, and a certain condition on
$\Gamma$ is satisfied (the ``semigroup condition''), then one can
use the associated splice diagram $\Delta$ to write down explicit
equations of a
complete intersection singularity, whose link is what we started
with.  This works as follows: first, for each pair of distinct
vertices $v,v'$ of $\Delta$, define the linking number $\ell_{vv'}$ to
be the product of the weights adjacent to, but not on, the shortest
path from $v$ to $v'$ (including weights around each vertex).  To
each leaf $w$, assign a variable $z_{w}$.  To each node $v$, assign
a weight $\ell_{vw}$ to the variable $z_{w}$, and assign a weight to
$v$ itself equal to the product of all the weights on the edges
adjacent to $v$.  Next, for each node $v$ and adjacent edge $e$, choose
 \emph{if possible}  a monomial $M_{ve}$ in the outer
variables whose total weight  is the weight of the node $v$.  If the
node $v$ has valency $\delta_{v}$, choose $\delta_{v}-2$ equations
 by equating to $0$
some $\C$-linear combinations of these monomials:
$$\sum_e a_{ie}M_{ve}=0,\quad i=1,\dots,\delta_v-2.$$
Repeating for all nodes, we get a total of $t-2$ equations.  We
require the coefficients $a_{ie}$ of the equations be ``generic''
in the precise sense that all maximal minors of the
$(\delta_{v}-2)\times \delta_{v}$ matrix $(a_{ie})$ have full
rank.

\begin{example}
  For the $\Delta$ of the example above,
we associate variables
  $z_1,\dots,z_4$ to the leaves as follows:
  $$
  \splicediag{6}{10}{\\
    &&z_1&\Circ&&&&&&&&&\Circ&z_4\\
    \Delta\quad=&&&&&&\Circ\lineto[ulll]_(.25){2}\lineto[dlll]^(.25)3
    &&&\Circ\lineto[drrr]_(.25){5}\lineto[urrr]^(.25)2
    \lineto[lll]_(.2){11}_(.8){7}\\
    &&z_2&\Circ&&&&&&&&&\Circ&z_3 }\qquad$$
  At the left node, the weights of the variables turn out to be (in
  order) 21, 14, 12, 30, and the total weight at the node is 42;
  so possible monomials
  for the left node are $z_1^2$, $z_2^3$, and $z_3z_4$. The
  monomials for the right node are $z_3^5$, $z_4^2$, and $z_1z_2^4$
  or $z_1^3z_2$.  Thus the system of equations might
  be
$$
\begin{array}{r}
  z_1^2+z_2^3+z_3z_4=0\,,\\
  z_3^5+z_4^2+z_1z_2^4=0\,.
\end{array}
$$

\end{example}

The ``semigroup condition'' on $\Delta$ (or $\Gamma$) is exactly the ability to
write down appropriate monomials at every node in every direction.

\begin{theorem}  Let $\Delta$ be a splice diagram corresponding to
a $\mathbb Z HS$ singularity link.  Suppose $\Delta$ satisfies the
semigroup conditions.  Then the splice equations above describe a
complete intersection singularity whose link is the $\mathbb Z HS$
associated to $\Delta$.  Further, each of the $t$ coordinates,
when set equal to 0, cuts out the knot on the link corresponding
to the end (``leaf'') of the diagram.
\end{theorem}

So, we can summarize by saying that for all these topologies,
there exists a very special kind of complete intersection
singularity with the given link.  This notion is a generalization
of Brieskorn complete intersection $V(p_{1}, \cdots ,p_{t})$
already discussed.  But in fact we can (and therefore should)
generalize the above construction slightly by allowing one to add
higher weight terms to each equation at each node.  We then arrive
at the notion of a \emph{complete intersection of splice type}, or
CIST.  The only proof we know of this Theorem is as a special case
of a much more general result, Theorem \ref{main} below.

We also note that each node in the splice diagram corresponds to a
valuation in the local ring of the CIST.  For, the nodes give weights
to the variables, and the nature of the defining equations means that
the associated graded ring is an integral domain (follows from
\cite{nw1}, Theorem 2.6).  This parallels the role of the valuation
for an irreducible curve, and the weight filtration for a weighted
homogeneous singularity.

The following (Example 3 of \cite{nw2}) shows that while the
embedding dimension of a CIST is at most the number of ends of
$\Delta$, it could be considerably smaller.

\begin{example}
 Let $\Delta$ be the splice diagram:
$$\splicediag{15}{20}{
\lefttag\Circ y {6pt}\lineto[drr]^(.75){q}&&&&&&&&\righttag\Circ
z {6pt}\lineto[dll]_(.75){p'}\\
&&\Circ\lineto[rr]^(.25){p''q'}^(.75){p}&&\Circ
\lineto[rr]^(.25){p'}^(.75){pq''}\lineto[dd]^(.25){p''}^(.75){p'qr}&&\Circ\\
\lefttag\Circ x {6pt}\lineto[urr]_(.75)p&&&&&&&&\righttag\Circ
w {6pt}\lineto[ull]^(.75){q'}\\
&&&&\Circ\lineto[ddl]_(.25){q''}\lineto[ddr]^(.25){p''}\\ \\
&&&\lefttag\Circ v {6pt}&&\righttag\Circ u{6pt}
}
$$
The integers $p$, $q$, $p'$, $q'$, $p''$, $q''$, $r$ are $\ge2$
and must satisfy appropriate relative primeness conditions, as
well as edge inequalities $$q'>p'q,\quad q''>p''q',\quad
qr>pq''\,.$$ Associating variables $x,y,z,w,u,v$ to the leaves in
clockwise order starting from the left as shown, one may write
splice equations:
\begin{align*}
  x^p+y^q&=z\\
z^{p'}+w^{q'}&=u\\
u^{p''}+v^{q''}&=x^r\\
y+w&=v
\end{align*}
These define the hypersurface singularity given by
$$((x^p+y^q)^{p'}+w^{q'})^{p''}+(y+w)^{q''}=x^r\,.$$
\end{example}

Given our general earlier warnings about NSS's, we can ask which of
the analytic types $(Y,0)$ for a given $\mathbb Z HS$ topology are so represented.
Of course, $V$ must be Gorenstein if it is a CIST.  But a very natural point
of view
is that one should be considering the ``algebraic'' nature not just of
the link, but the $t$ isotopy classes of knots which are part of the
data of the link.  We mentioned above that for a CIST, these knots are
cut out by coordinate functions.  Algebraically, this says that in the analytic
local ring of the singular point, there are prime principal ideals
which give topologically the knots in question.  A converse statement
holds:

\begin{theorem}\label{zhs}  Let $(Y,0)$ be a NSS with $\mathbb Z$HS link.
Suppose each of the $t$ knots in the link is represented by the
vanishing of some function in the local ring.  Then $Y$ is a CIST; in
particular, the link satisfies the semigroup condition.
\end{theorem}

The method of the proof is as follows: choose an irreducible curve
in the local ring cut out by the function corresponding to one of
the knots.  The other functions have a known order of vanishing
along the normalization of this curve, hence contribute to the
value semigroup.  This subsemigroup, read off from the splice
diagram, is shown to satisfy a certain inequality between its own
$\delta$ invariant and the ``Milnor number'' of the curve itself.
Applying now basic results of Buchweitz-Greuel
\cite{buchweitz-greuel}, one proves that the subsemigroup is the
full value semigroup of the curve, and these functions generate
the maximal ideal.  (Note that we did not need to assume the
Gorenstein property at the beginning.)

These results clarify greatly the possible nice geometries for a NSS
with given $\mathbb Z HS$ link; we even know how to write down explicit
equations for singularities whose link satisfies the semigroup
condition.  On the other hand, one should keep in mind the following
examples and open questions:

\begin{enumerate}
    \item Does every complete intersection singularity with $\mathbb
    Z HS$ link satisfy the semigroup condition?  Is every one
    a CIST? (This is related to a question about the Casson
    invariant.)
    \item Does every CIST as above satisfy the Neumann-Wahl
    Casson Invariant Conjecture?
    \item There exist Gorenstein singularities with $\mathbb Z HS$
    link which do not satisfy the semigroup condition (\cite{nemethi et al}, 4.5).
    \item There exists a Gorenstein singularity, not a complete
    intersection, whose link is $\Sigma(2,13,31)$ (\cite{nemethi et al}, 4.6).
    \end{enumerate}

    The second item above is difficult simply because one knows of
    no good inductive way to compute the geometric genus of a
    complete intersection of splice type, even though the
    equations are quite explicit.

   \section{generalized splice diagrams and CIST's}

    The results of the preceding section say a great deal about
    possible equations for a wide class of integral homology sphere
    links.  But from one point of view, such links are not so
    common.  Specifically, all rational singularities and nearly all
    minimally elliptic singularities have rational homology sphere
    link, but only a few have $\mathbb Z HS$ link.  In the rational
    case, the only non-trivial example is the $E_{8}$ singularity
    $V(2,3,5)$ (this is a famous theorem of Brieskorn, usually stated
    in terms of trivial local divisor class group).  Among minimally
    elliptics, one has $V(2,3,7)$ and $V(2,3,11)$ and their positive
    weight deformations.  It is thus natural to try to extend the previous discussion
    of CIST's to say something about a NSS with $\mathbb Q HS$ link.

   Given $(Y,0)$ with $\mathbb Q HS$ link and diagram
   $\Gamma$, one would like to get hold of the UAC $(X,0)\rightarrow
   (Y,0)$.  If $X$ is to be a complete intersection of a special type,
   one should try a generalization of the CIST's of the last section.
   In that case, one started with a splice diagram satisfying certain
   rules, and asked whether a certain ``semigroup'' condition was
   satisfied; then, one could write down equations of a
   complete intersection surface singularity, whose topology was what
   one wanted.

   Let us consider a more general splice diagram $\Delta$, where the weights
   around a node are positive, but are no longer required to be pairwise
   coprime.  (For technical reasons, one should also allow 1 to be a weight
   on an edge leading to a leaf.)  Then exactly as before, one can
   associate a variable to each of the $t$ ends; for each node,
   assign weights to the variables and the node; choose (again, \emph{if
   possible}) for each node and adjacent edge a monomial in the outer
   variables whose weight is that of the node; for each node, take generic (in a very
   specific sense) linear combinations of these monomials, giving
   weighted homogeneous polynomials for the node's weights; to each
   such polynomial, add terms of higher weight; set all these
   polynomials equal to 0.  In  other words, if the ``semigroup''
   condition is satisfied for a general splice diagram $\Delta$, one
   can as before produce subschemes $X(\Delta)$ of $\mathbb C ^{t}$.
   Then a major result of \cite{nw1} is

   \begin{theorem}  Suppose $\Delta$ is a generalized splice diagram
   satisfying the semigroup condition.  Then $X(\Delta)$ has an
   isolated local complete intersection surface singularity.
   \end{theorem}

 These singularities, which we still call CIST's (complete
 intersections of splice type), are the desired generalizations of
 Brieskorn complete intersections.  In fact, when the splice diagram
 has one node, an $X(\Delta)$ is exactly a BCI but with higher weight
 terms possibly added to each equation.  On the other hand, it is far
 from obvious how to prove that $X(\Delta)$ has an isolated
 singularity, at least with the very specific genericity condition we
 impose on the coefficient matrices.  This is accomplished in
 \cite{nw1} by an induction on the number of nodes.

 Starting with a graph $\Gamma$ representing a $\mathbb Q HS$, one can
 produce formally a generalized splice diagram $\Delta$ by the same
 procedure as in the $\mathbb Z HS$ case: collapse all vertices of
 valency 2, place certain subdeterminants as weights along every edge emanating
 from a node.  The differences now are that the weights around a node
 need not be pairwise relatively prime if one did not start with a
 $\mathbb Z HS$; and, the constructed splice diagram has less obvious
 topological interpretation than in the earlier case.
 Further, it is easy to see that different $\Gamma$'s can give rise to
 the same splice diagrams (this happens already in the weighted homogeneous case).
 However, we do have an unpublished result which indicates one is on
 the right track (and compares with Neumann's Theorem \ref{UAC}).

 \begin{theorem} Suppose two $\mathbb Q HS$ links give rise to the
 same splice diagram.  Then these two links have diffeomorphic
 universal abelian covers.
 \end{theorem}

 Returning to our original singularity $(Y,0)$, we have produced from
 the graph $\Gamma$ a splice diagram $\Delta$, and from it a class of
 isolated complete intersection surface singularities.  We hope one
 of these could be the UAC of $Y$.  But we need to bring into play the first
 homology group of the link; this ``discriminant group'' is computed
 from the intersection matrix $(E_{i}\cdot E_{j})$, and will be denoted
 $D(\Gamma)$.  If $\mathbb E$ denotes the free abelian group generated
 by the exceptional divisors on $\tilde{Y}$, then the intersection
 pairing gives an injective map of free $\mathbb Z$-modules $$\mathbb
 E \rightarrow \mathbb E ^{*}=\text{Hom}(\mathbb E ,\mathbb Z),$$
 whose cokernel is the discriminant group.

 \begin{proposition} Let $\{E_{i}\}$ be the exceptional curves, and
 $\{e_{i}\}\subset \mathbb E ^{*}$ be the dual basis for the
 intersection pairing, i.e. $$e_{i}\cdot E_{j}=\delta_{ij}.$$
 Then a faithful diagonal representation of the
 discriminant group on $\mathbb C ^{t}$ (the vector space with basis
 the ends of the graph) is constructed as follows:  $e\in \mathbb E ^{*}$
 acts on the coordinate $z_{i}$ by multiplication by the root of
 unity exp($2\pi i (e\cdot e_{i})$), where $e_{i}$ is dual basis
 element corresponding to the end $E_{i}$.
 \end{proposition}

 In other words, one has a natural
 representation of the discriminant group on the polynomial ring in
 $t$ variables, the ring from which the CIST's can be defined.  So,
 we ought to look for \emph{some} CIST on which the discriminant group
 acts equivariantly, i.e., for which every term of each defining
 equation transforms by the same character of the group.
 The semigroup condition guaranteed the existence of at least
 one ``admissible''
 monomial for each node and adjacent edge; we need to be able to find one that
 transforms correctly.   This translates easily into a
 condition on the original graph $\Gamma$, which we call the
 \emph{congruence condition}.  We can state the main theorem of
 \cite{nw1}.

 \begin{theorem}\label{main}  Let $\Gamma$ be a graph of a $\mathbb Q HS$ link
     satisfying the semigroup
 and congruence conditions, with associated splice diagram $\Delta$.
 Let $X(\Delta)$ be a complete intersection of splice type on which
 the discriminant group $D(\Gamma)$ acts equivariantly.  Then
 \begin{enumerate}
     \item $D(\Gamma)$ acts freely on $X(\Delta)$ off the singular
     point at the origin
     \item $Y\equiv X(\Delta)/D(\Gamma)$ is a germ of a NSS,
     whose resolution graph is $\Gamma$.
     \item $X(\Delta)\rightarrow Y$ is the universal abelian covering.
     \end{enumerate}
     \end{theorem}

 The bottom line is that if we are given a graph $\Gamma$
     satisfying the semigroup and congruence conditions, then we can
     ``explicitly'' write down the equations of a singularity with that
     link, in much the same way as discussed in the weighted homogeneous case
     following Citex.x.  That is, we can write
     explicit equations of a complete intersection singularity (the
     UAC), and an explicit diagonal action of the discriminant group
     on that singularity.  To see the actual equations of the
     desired singularity, one needs to do (perhaps very
     complicated) calculation of monomial invariants for the group action,
     find generators for the ideal of relations, and then deduce relations for
     these invariants which come from the splice equations.
     The easy case of Example \ref{easy} on the $E_{7}$ singularity as a
     quotient of the $E_{6}$ gives the general
     idea.  We have already mentioned that in the weighted homogeneous
     cases, there are faster ways to get equations than the UAC method
     of Neumann's Theorem \ref{UAC}.

      A singularity $Y$ arising as in the Theorem is called a \emph{splice
 quotient}.   A natural question is to ask which singularities are of this
 type.  We know that weighted homogeneous singularities
 with $\mathbb Q HS$ link are splice quotients.  Theorem \ref{zhs} gives an
 analytic necessary and sufficient condition for a singularity with
 $\mathbb Z HS$ link to be a splice quotient.  On the other hand,
 an ``equisingular'' deformation of a splice quotient need not be
 of that type; even if the geometric genera for the
 singularities in a family are constant, the same need not be true
 for the geometric genera of the UAC's.  An example of this
 phenomenon is found in \cite{nemethi et al}.

 Nonetheless, we conjectured about 7 years ago that rational and
 $\mathbb Q HS$ link
 minimally elliptic singularities are all splice quotients.  (By the
 time \cite{neumann-wahl02} was written, we had intemperately generalized the
 conjecture to a point where it could not be correct, via 
 \cite{nemethi et al}.)  The first
 non-trivial case was verified in \cite{neumann-wahl03} for the ``quotient cusps,''
 a class of log-canonical (and taut) rational singularities, whose resolution
 dual graph has 2 nodes:
 $$
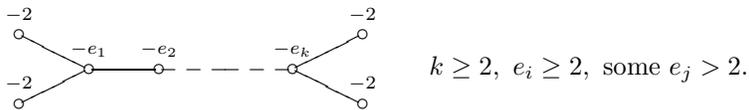
\xymatrix@R=6pt@C=24pt@M=0pt@W=0pt@H=0pt{
  \overtag{\Circ}{-2}{8pt}\lineto[dr] && &&&
  \overtag{\Circ}{-2}{8pt}\lineto[dl]\\
  &\overtag{\Circ}{-e_1}{8pt}\lineto[r]
  &\overtag{\Circ}{-e_2}{8pt}\dashto[r]&\dashto[r]&
  \overtag{\Circ}{-e_k}{8pt}&&k\ge2,~~e_i\ge2, \text{ some }e_j>2.\\
  \overtag{\Circ}{-2}{8pt}\lineto[ur] && &&&
  \overtag{\Circ}{-2}{8pt}\lineto[ul]}$$
  Explicit equations for the UAC (which is a ``cusp'' singularity)
  and the action of the discriminant group are given in Section 5 of
  that paper.

 The motivation for the general conjecture was not only
 the beauty of such a result, but because the rational and $\mathbb Q
 HS$ minimally elliptic
 singularities possess an important property
 of splice quotients analogous to that mentioned in Theorem x.x
 above.  Recall that an \emph{end-curve} on the minimal good
 resolution $\tilde{Y}$ is a rational curve that has just one
 intersection point with the rest of the exceptional divisor (so that
 it corresponds to a leaf in the splice diagram).  The following is observed
 along the way to proving the Main Theorem above.

 \begin{proposition}\label{endc} Let $(Y,0)$ be a splice quotient.  Then for every
 end curve $E_{i}$ on $\tilde{Y}$, there is a function
 $y_{i}:Y\rightarrow \mathbb C$ such that the proper transform on
 $\tilde{Y}$ of its zero-locus consists of one smooth irreducible
 curve $C_{i}$, which intersects $E_{i}$ transversally at one point
 and intersects no other exceptional curve.

 \end{proposition}

 Another way to state this property is that for every end-curve,
 there is a prime ideal in the analytic local ring of $Y$
 whose $n_{i}$-th symbolic
 power is a principal ideal $(y_{i})$, where $y_{i}$ has the vanishing
 properties described above (i.e., its proper transform is
 $n_{i}C_{i}$.)  Note that this integer $n_{i}$ is the order of the
 image of the dual basis element $e_{i}$ in the divisor class group.

 Now, it is well-known that rational singularities have the
 ``end-curves property'' described in the Proposition; the same is
 true for $\mathbb Q HS$ link minimally elliptics (\cite{reid}, p. 112).
 So, in an attempt to generalize Theorem \ref{zhs} we have made the
 following

 \begin{ecconjecture*} Suppose $(Y,0)$ is a NSS with $\mathbb Q HS$
 link.  Suppose to every end-curve on the minimal good resolution
 there exists a function as in Proposition \ref{endc}.  Then $Y$ is a splice
 quotient.
 \end{ecconjecture*}

 Note that the assumptions about the end-curves are supposed to imply
 the semigroup and congruence conditions on the graph, as well as
 the $\mathbb Q$-Gorensteinness of the singularity.

 This Conjecture is still open; but as we shall see, T. Okuma has
 recently proved that rational and minimally elliptic singularities
 are splice quotients.
 
 (Remark added, September 2005: the author and W. Neumann have 
 announced a proof of this Conjecture.)

 \section{Okuma's theorem and further questions}

 Once we know that a graph $\Gamma$ of a rational or $\mathbb Q HS$
 minimally elliptic singularity satisfies the semigroup and
 congruence conditions, then it follows from Theorem \ref{main} that there is at least
 one such singularity which is a splice quotient.  For $\Gamma$ with
 at most two nodes, these two conditions may be checked directly
 (\cite{nw1},Section 11).   But T. Okuma has proved in general

 \begin{proposition}\cite{okuma2}  The graph of a rational or
 $\mathbb Q HS$ minimally
  elliptic singularity satisfies the semigroup and congruence
  conditions.
  \end{proposition}

  Okuma's method is to give a condition on $\Gamma$ that turns out to
  be equivalent to the semigroup and congruence conditions, and then
  to deduce this graph-theoretic property from a well-known stronger
  property for such rational or $\mathbb Q HS$ minimally elliptic
  graphs.  Our own version of his result is found in \cite{nw1},
  Section 13.

  To get the strongest result, Okuma uses a precise
  description of the UAC of a NSS $Y$ with $\mathbb Q HS$ link.
  Considering the MGR $\tilde{Y}\rightarrow Y$, he constructs a
  fairly explicit sheaf of algebras on $\tilde{Y}$ whose Spec is a
  partial resolution of the UAC, with only cyclic quotient singularities \cite{okuma1}.
    This is similar to the Esnault-Viehweg method for
  constructing cyclic coverings branched along normal crossings
  divisors.  Using the preceding proposition, and the existence of
  appropriate end-curves, Okuma proves our old Conjecture about the UAC.

  \begin{theorem}(\cite{okuma2})  Every rational or $\mathbb Q HS$ minimally
  elliptic singularity is a splice quotient.  In particular, one may
  write down explicit equations for it.
  \end{theorem}

  We note that Okuma's original preprint does not specifically assert
  this Theorem in its full strength; one can find an explanation of
  why he has in fact obtained this result in \cite{nw1}, Section 13.

  At this point, we now know many examples of singularities with $\mathbb Q HS$ links which
  are splice quotients---especially, rational, minimally elliptic,
  and weighted homogeneous.  But there are many examples, even of
  hypersurface singularities, which could not be splice quotients.
   The next challenge is to try to understand better what is going
   on in the other cases---is there a nice theorem out there?

   An interesting place to start is with some of the examples in
   \cite{nemethi et al}.  For instance, we have found a hypersurface
   singularity which does not satisfy the semigroup condition; but
   nonetheless, the UAC is a complete intersection of splice type!
    These and related issues are currently being looked into.



\begin{thebibliography}{99}

%
\bibitem{buchweitz-greuel} Ragnar-Olaf Buchweitz, Gert-Martin Greuel,
The Milnor number and deformations of complex curve singularities.
Invent. Math. {\bf58} (1980), 241--281.
%
%
%
%


\bibitem{eisenbud-neumann} D. Eisenbud and W.D. Neumann,
  \emph{Three-dimensional link theory and invariants of plane curve
singularities.}  Ann. Math. Stud. {\bf 110}, Princeton.  Princeton
  Univ. Press (1985).

%
%
%
%
%
%
%
%
\bibitem{ds} T. de Jong and D. van Straten,
On the deformation theory of rational surface singularities
with reduced fundamental cycle. J. Algebraic Geom. {\bf 3} (1994), 117--172.


\bibitem{laufert} H.B. Laufer, Taut two-dimensional singularities.
 Math. Ann. {\bf 205} (1973), 131-164.

\bibitem{laufermin} 
\bysame,
On minimally elliptic singularities. Amer. J. Math.
{\bf99} (1977),  1257--1295.

\bibitem{laufermu} 
\bysame,
On $\mu$ for surface singularities,
 {\it Several complex variables},
 Proc. Symp. Pure Math. {\bf 30} (Amer. Math. Soc. 1977), 45-49

\bibitem{milnor-book} John Milnor, \emph{Singular points of complex
    hypersurfaces}, Annals of Mathematics Studies, {\bf61} (Princeton
  University Press, Princeton, N.J.; University of Tokyo Press, Tokyo
  1968)

\bibitem{nemethi et al} I. Luengo-Velasco, A. Melle-Hern\'andez, A.
  N\'emethi, Links and analytic invariants of superisolated
  singularities. J. Algebraic Geom. {\bf 14} (2005), 543--566.


%
%
%
%
\bibitem{neumann81} W.D. Neumann,  
A calculus for plumbing applied to the topology of complex
  surface singularities and degenerating complex curves. Trans. Amer.
  Math. Soc. {\bf 268} (1981), 299--343.

\bibitem{neumann83} 
\bysame,
Abelian covers of quasihomogeneous
  surface singularities, {\it Singularities, Arcata 1981}, Proc. Symp.
  Pure Math. {\bf 40} (Amer. Math. Soc. 1983), 233--243.

\bibitem{n} 
\bysame, Topology of hypersurface singularities, {\it Erich
K\"ahler--Mathematische Werke}, R. Berndt and O. Riemenschneider,
eds. (Walter de Gruyter Verlag, 2003), 727--736.

\bibitem{neumann-wahl90} W.D. Neumann, J. Wahl, Casson invariant of links of
singularities. Comment. Math. Helv. {\bf 65} (1990), 58--78.


\bibitem{neumann-wahl02} 
\bysame,
Universal abelian covers of surface singularities, {\it
Trends on Singularities}, A. Libgober and M. Tibar, eds.
  (Birkh\"auser Verlag, 2002), 181--190.

\bibitem{neumann-wahl03} 
\bysame,
Universal abelian covers of quotient-cusps.
  Math. Ann. {\bf 326} (2003), 75--93.


\bibitem{nw1}
\bysame,
Complete intersection singularities of splice
 type as universal abelian covers. Geom. Topol. {\bf 9} (2005), 699-755.

\bibitem{nw2} 
\bysame,
Complex surface singularities with integral homology sphere
 links. Geom. Topol. {\bf 9} (2005), 757-811.


\bibitem{okuma1} T. Okuma, Universal abelian covers of rational
 surface singularities. J. London Math. Soc. (2) {\bf 70} (2004), 307-324.

\bibitem{okuma2} T. Okuma, Universal abelian covers of certain
 surface singularities, http://arxiv.org/math.AG/0503733

 \bibitem{O-W} P. Orlik and P. Wagreich, Isolated singularities of algebraic
 surfaces with $\mathbb C ^{*}$ action.  Ann. of Math. {\bf 93} (1971), 205--228.



\bibitem{pinkham} H. Pinkham, Normal surface singularities with
$\C ^{*}$ action.  Math. Ann. {\bf 227} (1977), 183--193.

\bibitem{reid} M. Reid, Chapters on Algebraic Surfaces, \emph{Complex
algebraic geometry (Park City UT, 1993)}, IAS/Park City Math. Ser.
3, Amer. Math. Soc. (Providence, RI, 1997), 3-159.

%
%


\bibitem{wahl} J. Wahl, Equations defining rational singularities.
Ann. Sci. \'Ecole Norm. Sup. {\bf 10} (1977), 231-264.


\bibitem{zariski} Oscar Zariski, \emph{Le probl\`eme des modules pour
   les branches planes.} (\'Ecole Polytechnique, Paris, 1973)

\end{thebibliography}
\end{document}